\documentclass[]{iac}

\newtagform{brackets}{(}{)}
\usetagform{brackets}

\usepackage{graphicx}
\usepackage{amsmath}
\usepackage{bm}
\usepackage[version=4]{mhchem}
\usepackage[per-mode = symbol]{siunitx}  
\DeclareSIUnit\permille{\text{\textperthousand}}

\usepackage[mode=buildnew]{standalone}
\usepackage{tikz}

\usepackage{hyperref} 
\hypersetup{
    colorlinks = true,  
	allcolors = black  
}

\RequirePackage{color}
\usepackage{xcolor}
\usepackage{soul}

\usepackage[super]{nth}
\usepackage{mathtools}

\usepackage{enumitem}

\usepackage{footnote}
\usepackage{threeparttable, tablefootnote}
\usepackage{multirow}
\usepackage{booktabs}

\usepackage{caption}
\usepackage{subcaption}

\allowdisplaybreaks

\newcommand{\real}{\mathbb{R}}

\newcommand{\ev}[1]{\mathbb{E}[#1]}

\newcommand{\pr}[1]{\text{Pr}\{#1\}}

\newcommand{\abs}[1]{\lvert#1\rvert}

\newcommand{\norm}[1]{\left\lVert#1\right\rVert}

\newcommand{\skewmat}[1]{\left[{#1}_\times\right]}

\def \StageThree {1}
\def \CoastingThree {2}
\def \StageFourOne {3}
\def \CoastingFour {4}
\def \StageFourTwo {5}
\def \Return {6}

\def \tf {t_f}

\begin{document}

\IACpaperyear{22}
\IACpapernumber{D2.IPB.12.x71527}
\IACconference{73}
\IAClocation{Paris, France}
\IACdate{18--22 September 2022}
\IACcopyrightB{2022}{Dr. Boris Benedikter}

\title{Stochastic Control of Launch Vehicle Upper Stage with Minimum-Variance Splash-Down}

\IACauthor{Boris~Benedikter$^{\text{a}\ast}$, }{$^\text{a}$\textit{Postdoctoral Researcher, Department of Mechanical and Aerospace Engineering, Sapienza University of Rome, Via Eudossiana 18, 00184, Rome, Italy}, \uline{boris.benedikter@uniroma1.it}}

\IACauthor{Alessandro~Zavoli$^{\text{b}}$, }{$^\text{b}$\textit{Assistant Professor, Department of Mechanical and Aerospace Engineering, Sapienza University of Rome, Via Eudossiana 18, 00184, Rome, Italy}, \uline{alessandro.zavoli@uniroma1.it}}

\IACauthor{Guido~Colasurdo$^{\text{c}}$, }{$^\text{c}$\textit{Full Professor, Department of Mechanical and Aerospace Engineering, Sapienza University of Rome, Via Eudossiana 18, 00184, Rome, Italy}, \uline{guido.colasurdo@uniroma1.it}}

\IACauthor{Simone~Pizzurro$^{\text{d}}$, }{$^\text{d}$\textit{Research Fellow, Launchers and Space Transportation Department, Italian Space Agency, Via del Politecnico snc, 00133, Rome, Italy}, \uline{simone.pizzurro@est.asi.it}}

\IACauthor{Enrico~Cavallini$^{\text{e}}$}{$^\text{e}$\textit{Head of Space Transportation Programs Office, Space Transportation, Space Infrastructures, and In-Orbit Servicing Department, Italian Space Agency, Via del Politecnico snc, 00133, Rome, Italy}, \uline{enrico.cavallini@asi.it}}

\abstract{%
This paper presents a novel  synthesis method for designing 
an optimal and robust guidance law for a non-throttleable upper stage of a launch vehicle, using a convex approach.
In the unperturbed scenario, a combination of lossless and successive convexification techniques is employed to formulate the guidance problem as a sequence of convex problems that yields the optimal trajectory, to be used as a reference for the design of a feedback controller,
with little computational effort.
Then, based on the reference state and control, a stochastic optimal control problem is defined to find a closed-loop control law that rejects random in-flight disturbance.
The control is parameterized as a multiplicative feedback law; thus, only the control direction is regulated, while the magnitude corresponds to the nominal one, enabling its use for solid rocket motors.
The objective of the optimization is to minimize the splash-down dispersion to ensure that the spent stage falls as close as possible to the nominal point.
Thanks to an original convexification strategy, the stochastic optimal control problem can be solved in polynomial time since it reduces to a semidefinite programming problem.
Numerical results assess the robustness of the stochastic controller and compare its performance with a model predictive control algorithm via extensive Monte Carlo campaigns.
}

\maketitle

\section{Introduction}

Launch vehicle dynamics are subject to significant uncertainties due to dispersions of the propulsion system performance, random variations of the local environment, and hard-to-model external perturbations.
Therefore, the design of a guidance and control algorithm that can robustly ensure operation under uncertainty is essential for the system's in-flight autonomy.
Besides safety requirements, the algorithm should also be optimal, e.g., in terms of propellant consumption, to maximize the mass-to-orbit, enhance the vehicle responsiveness in off-nominal conditions, and reduce the mission cost.

The optimal control of uncertain systems is traditionally an uneasy task, as a general approach to an uncertain optimal control problem (OCP) would involve solving a dynamic programming (DP) problem over arbitrary feedback control laws \cite{bellman1966dynamic}.
However, DP is computationally impractical for real-world applications due to the curse of dimensionality.
Differential dynamic programming (DDP) can yield an optimal control policy with lower computational complexity than DP by expanding the OCP around a reference solution \cite{jacobson1970differential},
but it requires additional effort to retrieve the analytical formulation of the second-order expansion of the objective function, large disturbances may deviate excessively the system from the reference solution, and the inclusion of nonlinear constraints is generally nontrivial \cite{ozaki2020tube}.


The interest in machine learning (ML) techniques grew significantly in the last few years for real-time guidance applications thanks to the low evaluation times and the high accuracy in function approximation of deep neural networks (DNNs) \cite{izzo2018survey}.
In particular, a DNN can be trained to imitate an expert behavior (hence, a robust control policy) through behavioral cloning \cite{wang2020asteroid,federici2021deep}.
Also, reinforcement learning was proposed to train DNNs via repeated interactions with the environment (i.e., simulations of the stochastic dynamics) to progressively learn a robust optimal control policy \cite{zavoli2021reinforcement, federici2021autonomous, federici2021reinforcement}.
However, the lack of any theoretical guarantees on the robustness of the policy and the extensive computational effort associated with the DNN training process, still limit the application of ML methods to relatively simple problems.

The potentially prohibitive computational cost of DP or DNN training can be bypassed if leveraging a model predictive control (MPC) scheme.
Indeed, MPC, rather than relying on a feedback policy computed pre-flight, recursively solves an OCP onboard to update the control signal based on real-time measurements.
The recursive update provides the algorithm with inherent robustness to in-flight disturbances and model errors even without necessarily accounting for uncertainty in the OCP.
The ability to simultaneously maximize system performance and enforce mission requirements made MPC extremely appealing for aerospace guidance and control applications \cite{eren2017model, benedikter2020autonomous, benedikter2021autonomous, iannelli2021attitude}.
However, the degree of robustness of an MPC algorithm depends on the update frequency, which is limited by the efficiency of the onboard hardware in solving the OCP.

In this respect, convexification techniques gained great popularity in the last few years, as they allow for solving originally nonconvex problems with highly efficient convex programming algorithms \cite{liu2017survey}.
For instance, \emph{lossless} convexification can be used to formulate problems with a quite common class of nonconvex control constraints as equivalent convex programs solvable in polynomial time \cite{accikmecse2011lossless}.
More general sources of nonconvexity can be tackled via \emph{successive} convexification, which comes with less sound theoretical guarantees, but proved to be successful in numerous aerospace problems, including 
spacecraft rendezvous~\cite{benedikter2019convexrendezvous},
rocket ascent \cite{benedikter2019convexascent,benedikter2021convex, benedikter2022convex} 
and landing \cite{szmuk2016successive,sagliano2018pseudospectral},
low-thrust transfers~\cite{wang2018optimization, wang2018minimum},
and atmospheric entry~\cite{wang2018autonomous, sagliano2018optimal}.
However, to ensure convergence, the parameters of the successive convexification algorithm must be carefully picked and some modifications to the formulation are often necessary to prevent undesired phenomena, such as artificial infeasibility or unboundedness, resulting in longer computational times and possibly limiting the update frequency (hence, the robustness) of an MPC architecture. 

To endow the optimal control policy with sound robustness guarantees, uncertainty must be included in the problem formulation.
In this respect, robust optimal control (ROC) methods account for uncertainty by modeling them as deterministic variables bounded within finite ranges.
In this way, the associated optimization problem is solvable with traditional programming algorithms.
However, this modeling strategy usually leverages either min-max OCP formulations \cite{bemporad1999robust} or constraint tightening \cite{chisci2001systems}, which often lead to overly conservative policies.
More recently, tube-based approaches have been proposed to design more performing policies than traditional ROC methods.
Specifically, by accounting for (bounded) uncertainties directly in the control algorithm,
tube-based approaches can find a control policy that encloses all future state realizations in a bounded region (referred to as a \emph{tube})
\cite{langson2004robust}.

On the other hand, stochastic optimal control (SOC) deals with unbounded uncertainties, which are the class of uncertainty that affect most real-world systems, by incorporating their probabilistic descriptions in the OCP \cite{mesbah2016stochastic}.
The unboundedness of the random variables requires formulating hard requirements on the state and controls as relaxed probabilistic constraints, called \emph{chance constraints}, which allow for a small probability of violation.
Chance-constrained optimization has been applied to several aerospace problems, 
such as UAV path planning with obstacle avoidance constraints \cite{blackmore2011chance}, 
model rocket ascent \cite{lew2019chance},
station-keeping in unstable orbits \cite{oguri2019convex}, 
and impulsive \cite{oguri2019risk} or finite-thrust \cite{oguri2019risk2} orbit transfers.
However, most works in the literature consider a simple open-loop control law that cannot control the time-evolution of the covariance, which increases uncontrolled over time due to diffusive stochastic perturbations.

To design performing solutions, the stochastic OCP should account for the presence of a feedback controller.
A simple approach consists in designing the feedback controller \emph{a priori} via traditional control theory approaches (e.g., Linear Quadratic Regulator, LQR) and then optimizing only the feedforward control signal in the stochastic framework \cite{lew2020chance}.
However, the quality of the solution depends on the considered controller, which, if not tuned properly, may still lead to over-conservative policies.
\emph{Covariance control}, originally proposed in the late 1980s for infinite-horizon problems \cite{hotz1987covariance}, gained renewed interest among the control community after being applied to finite-horizon problems as well \cite{chen2016optimal1,chen2016optimal2}.
The major appeal of this control synthesis method lies in the ability to systematically design, with low computational complexity, a closed-loop optimal control policy that drives a stochastically perturbed system from an initial probability distribution (in terms of mean state and covariance) to a desired one at a final time.
Indeed, for a linear system, the finite-horizon stochastic OCP can be cast as a deterministic convex optimization problem, thus solvable in polynomial time \cite{bakolas2016optimal1,bakolas2016optimal2}. 
Recently, thanks to convexification methods, covariance control has been applied to nonlinear systems as well, including low-thrust interplanetary trajectory design problem \cite{ridderhof2020chance,benedikter2022covariance}, powered descent and landing \cite{ridderhof2021minimum,benedikter2022landing}, and planetary entry \cite{ridderhof2022stochastic}.

In most works on covariance control, the control law is parameterized as the sum of a feedforward term and a linear feedback one.
Thus, when operating in uncertain environments, the system adjusts the control magnitude and direction to compensate for deviations from the nominal path.
However, such a control policy can only be implemented by a system that can regulate the magnitude of the control action.
Instead, non-throttleable systems must rely on policies that regulate exclusively the direction of the control vector, without straying from the nominal control magnitude profile.

This paper presents an algorithm to design a nominal optimal trajectory for a non-throttleable upper stage of a launch vehicle and a robust controller that compensates for exogenous disturbances.
The considered case study involves the flight of VEGA's third stage, which features a solid rocket motor (SRM), thus, its thrust level cannot be regulated in flight.
The system must be actively controlled to ensure that the spent stage falls in an uninhabited area;
this requirement is generally represented as a \emph{splash-down constraint}.
%

A multiplicative feedback control law that governs only the thrust vector direction is considered so that the closed-loop control norm is independent of the disturbances.
Due to the nonlinearity arising from the coupling of the feedback gains and the covariance matrix in the dynamics, the associated optimization problem is not convex.
Thus, to reduce the computational complexity of the solution process, the problem is solved in two steps.
First, the optimal trajectory and the associated  feedforward control signal are computed in the deterministic scenario, that is, neglecting uncertainties, and serve as a reference for the second step.
Then, a covariance control problem, based on the reference trajectory and feedforward control, is solved to find the multiplicative feedback gains that minimize the dispersion of the spent stage splash-down point.
In either step, the OCP can be solved with low computational complexity thanks to the introduction of convenient changes of variables, lossless constraint relaxations, and successive linearization.

The paper is organized as follows.
Section~\ref{sec:problem} describes the upper stage guidance and control problem, detailing the phases, the dynamical model, and the mission requirements of the corresponding OCP.
Section~\ref{sec:cvx} outlines the convexification strategy used to cast the nonlinear OCP as a sequence of convex problems that rapidly yields the nominal trajectory and control law in the unperturbed scenario.
Sections~\ref{sec:mpc} and \ref{sec:covariance_control} outline two different strategies to design a robust control framework for the upper stage flight.
In Section~\ref{sec:mpc}, 
an MPC algorithm that recursively solves the sequential convex optimization problem is proposed as a means to provide robustness to model mismatches and random disturbances.
Instead, Section~\ref{sec:covariance_control} outlines an optimal covariance control problem that yields a multiplicative feedback controller robust to the considered characterization of the in-flight perturbations.
Section~\ref{sec:results} presents numerical results to assess the robustness and performance of the proposed algorithms through extensive Monte Carlo campaigns and compares the MPC algorithm with the stochastic control approach.
A conclusion section ends the manuscript.


\section{Problem Statement}
\label{sec:problem}

In this section, the upper stage optimal control problem is presented.
First, the phase sequence of the OCP is outlined;
then, the considered dynamical model is presented;
and finally, the objective and all mission requirements of the OCP are discussed. 

\subsection{Phase Sequence}


The problem under investigation corresponds to the final portion of the VEGA launch vehicle ascent flight.
In particular, the operation of the first two stages is neglected and the goal of the optimization is to design an optimal control policy for the third stage flight.
To evaluate and maximize the vehicle performance in terms of carrying capacity the OCP horizon must extend until payload release.

\begin{figure}
    \centering
    \includegraphics[width=\columnwidth]{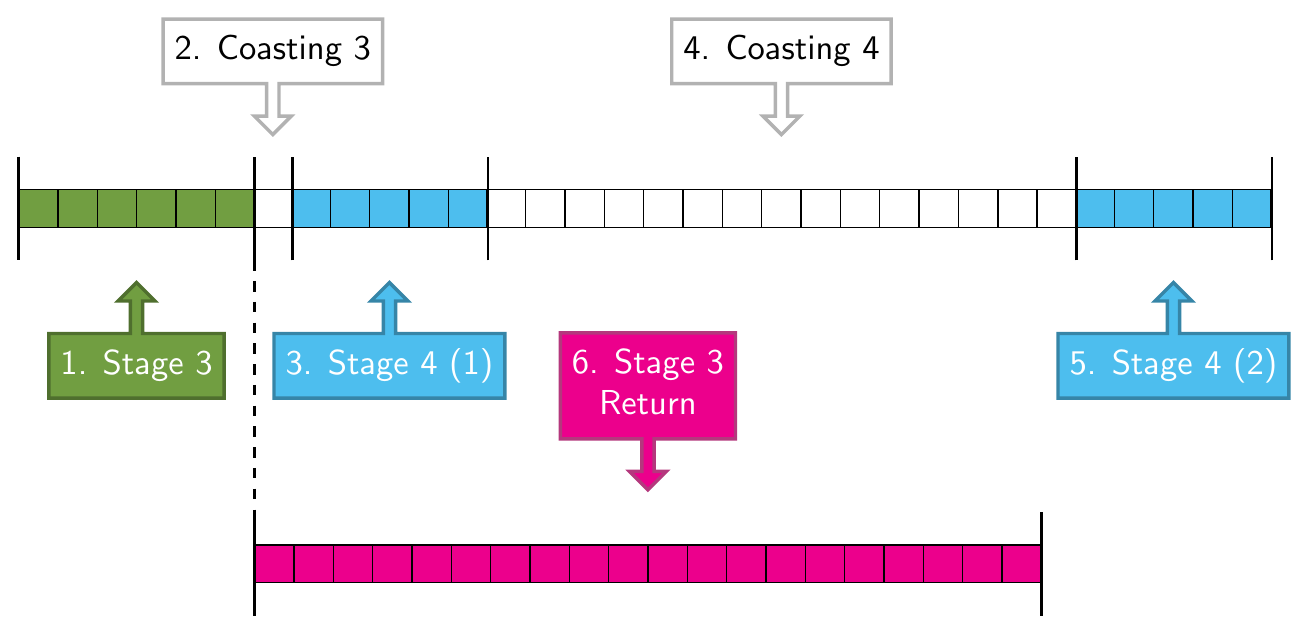}
    \caption{Phases of the optimal control problem.}
    \label{fig:phases}
\end{figure}

Figure~\ref{fig:phases} illustrates schematically the flight phases of the VEGA launch vehicle from the third stage ignition until payload injection into the target orbit.
The phases are numbered from \StageThree{} to \StageFourTwo{} in chronological order, and Phase~\Return{} simulates the uncontrolled return of the third stage after burnout.
The third stage flight corresponds to Phase~\StageThree{} and terminates at the burnout of the SRM.
A short coasting (Phase~\CoastingThree{}) of fixed duration follows the third stage separation.
Then, the fourth and final stage must complete the orbit insertion maneuver, whose optimal phase sequence is similar to a Hohmann maneuver, i.e., composed of two propelled arcs (Phases~\StageFourOne{} and \StageFourTwo{}) separated by a quite long coasting (Phase \CoastingFour{}).
The time-lengths of Phases~\StageFourOne{}--\StageFourTwo{} must be optimized.
The return (Phase~\Return{}) is included in the OCP to predict and constrain the splash-down point of the spent third stage.
Also the duration of Phase~\Return{} must be determined by the optimization process.



Hereinafter, let \smash{$t_0^{(i)}$} and \smash{$t_f^{(i)}$} denote the initial and final time of the $i$-th phase and $\smash{\Delta t^{(i)}}$ its time-length.
For the sake of simplicity, if no phase superscript is specified, then $t_0$ and $t_f$ denote the ignition time of the third stage $\smash{t_0^{(\StageThree{})}}$ and the fourth stage burnout $\smash{t_f^{(\StageFourTwo{})}}$, respectively.
Likewise, let $t_R$ denote the splash-down time of the spent stage \smash{$t_f^{(\Return{})}$}.

\subsection{System Dynamics}

The vehicle is modeled as a point mass subject to a 3-DoF translational motion. 
The rotational dynamics are neglected and the rocket axis is assumed to be aligned with the thrust direction at any time. 
The state $\bm{x}$ of the rocket is described by its position $\bm{r}$, velocity $\bm{v}$, and mass $m$, as $\bm{x} = [ x \; y \; z \; v_x \; v_y \; v_z \; m ]$.
The position and velocity are expressed in Cartesian Earth-Centered Inertial (ECI) coordinates. 
The $x$ and $y$ axes lie in the equatorial plane, and the $z$ axis is aligned with the Earth's angular velocity $\bm{\omega}_E$, forming a right-hand frame.
The main advantage of using this set of state variables over, for example, a spherical coordinate system \cite{benedikter2021convex}, is the fact that any trajectory can be studied, even missions to polar or high inclination orbits, which are the typical targets for VEGA, without running into singularities.

The only external forces assumed to act on the launch vehicle are gravity, the aerodynamic drag, and the motor thrust.
A simple inverse-square gravitational model is assumed, thus the gravity acceleration is $\bm{g} = - \mu \bm{r} / r^3 $,
with $\mu$ denoting the Earth's gravitational parameter.
The drag force is modeled as
$\bm{D} = - 0.5 \rho {v}_{\text{rel}}^2 S_{\text{ref}} C_D \bm{\hat{v}}_{\text{rel}}$,
where 
$\rho$ is the atmospheric density, 
$S$ is the reference surface, 
$C_D$ is the drag coefficient, assumed to be constant, 
and $\bm{v}_{\text{rel}} = \bm{v} - \bm{\omega}_E \times \bm{r}$ is the relative-to-atmosphere velocity.
Lift is neglected, as it is much smaller than drag for launch vehicles.

VEGA's third stage is an SRM, so the thrust history $T(t)$ and the associated mass flow rate $\dot{m}_e(t)$ depend on the geometry of the solid propellant.
Although upper stages operate at high altitudes, the net thrust is slightly different from the vacuum one due to the atmospheric pressure $p$ contribution, thus we consider a vacuum thrust law $T_{\text{vac}}(t)$ and a corresponding mass flow rate law $\dot{m}_e(t)$.
The net thrust is $T = T_{\text{vac}}(t) - p A_e$, where $A_e$ is the nozzle exit area.

Instead, the last stage is equipped with a small liquid rocket engine, AVUM; hence, its thrust and mass flow rate are assumed to be constant and equal to the maximum attainable values.
Also, AVUM can be cut off and re-ignited; so, its overall burn time can be minimized by the optimization process to save as much propellant as possible.

While the thrust magnitude is prescribed, the thrust direction $\bm{\hat{T}}$ is the control variable.
Its components are expressed in ECI coordinates and, being a unit vector, must satisfy the following identity
\begin{equation}
    \hat{T}_x^2 + \hat{T}_y^2 + \hat{T}_z^2 = 1
	\label{eq:thrust_direction_equality_path_con}
\end{equation}

The resulting equations of motion $\bm{\dot{x}} = \bm{f}(\bm{x}, \bm{u}, t)$ are
\begin{align}
    \dot{\bm{r}} &= \bm{v} \label{eq:original_ODE_r} \\
    \dot{\bm{v}} &= -\frac{\mu}{r^3} \bm{r} + \frac{T}{m} \bm{\hat{T}} - \frac{D}{m} \bm{\hat{v}}_{\text{rel}} \label{eq:original_ODE_v} \\
    \dot{m} &= -\dot{m}_e \label{eq:original_ODE_m}
\end{align}

\subsection{Mission Requirements and Objective}

The state at the third stage ignition is supposed to be completely assigned, that is
\begin{equation}
    \bm{x}(t_0) = \tilde{\bm{x}}_0
    \label{eq:x0}
\end{equation}
At the burnout of the fourth stage, the payload must be released into the target orbit.
In this work, we suppose that the target orbit is a circular orbit of semi-major axis $a_{\text{des}}$ and inclination $i_{\text{des}}$.
Thus, the final conditions to impose at the end of Phase \StageFourTwo{} are 
\begin{align}
    x(\tf)^2 + y(\tf)^2 + z(\tf)^2 &= a_{\text{des}}^2 \label{eq:final_radius_nonlinear} \\
    v_{x}(\tf)^2 + v_{y}(\tf)^2 + v_{z}(\tf)^2 &= v_{\text{des}}^2 \label{eq:final_velocity_nonlinear} \\
    \bm{r}(\tf) \cdot \bm{v}(\tf) &= 0 \label{eq:final_radial_velocity} \\
    x(\tf) v_y(\tf) - y(\tf) v_x(\tf) &= h_{z, \text{des}}
    \label{eq:final_ainc_circular}
\end{align}
Equation~\eqref{eq:final_radius_nonlinear} prescribes the final radius, while Eq.~\eqref{eq:final_velocity_nonlinear} with $v_{\text{des}} = \sqrt{\mu / a_{\text{des}}}$ and Eq.~\eqref{eq:final_radial_velocity} are sufficient to impose the null eccentricity constraint.
As for the inclination, $i_{\text{des}}$ is imposed by prescribing the $z$-component of the angular momentum vector $h_{z, \text{des}} = \cos i_{\text{des}} \sqrt{\mu a_{\text{des}}}$.

After burnout, the third stage features a high, but suborbital, velocity, so it falls back to Earth.
To prevent it from crashing into an inhabited area, Phase~\Return{} simulates its return and the splash-down point can be constrained by imposing the following conditions
\begin{align}
    x(t_R)^2 + y(t_R)^2 + z(t_R)^2 &= R_E^2 \label{eq:final_r_reentry_nonlinear} \\
    z(t_R) &= z_{R, \text{des}} \label{eq:final_LAT_reentry}    
\end{align}
with $z_{R, \text{des}} = R_E \sin\varphi_{R, \text{des}}$.
Equation~\eqref{eq:final_r_reentry_nonlinear} constrains the final radius to be equal to the Earth radius $R_E$, while Eq.~\eqref{eq:final_LAT_reentry} constrains the latitude of the splash-down point to a desired value $\varphi_{R, \text{des}}$.
For missions toward polar or quasi-polar orbits (e.g., Sun-synchronous orbits), which are the scope of the present paper, constraining the latitude is equivalent to constraining the distance from the launch site, which is the only quantity that can be actively controlled during the third stage flight.
Indeed, constraining the longitude to an arbitrary value is impractical, as it would be very expensive because it requires a maneuver to change the plane of the ascent trajectory.

The multi-phase structure of the problem can account for mass discontinuities, which take place at the third stage separation.
In particular, all state variables are continuous across phases, while the following linkage conditions must be imposed 
\begin{equation}
    m(t_0^{(\CoastingThree{})}) = m(t_f^{(\StageThree{})}) - m_{\text{dry}, 3} \label{eq:mass_lkg_con_3}
\end{equation}
As for the return phase, the initial position and velocity are equal to the burnout conditions of the third stage, while its mass corresponds to the dry mass of the stage.
Thus,
\begin{align}
    \bm{r}(t_0^{(\Return{})}) &= \bm{r}(t_f^{(\StageThree{})}) \label{eq:return_lkg_position} \\ 
    \bm{v}(t_0^{(\Return{})}) &= \bm{v}(t_f^{(\StageThree{})}) \label{eq:return_lkg_velocity} \\ 
    m(t_0^{(\Return{})}) &= m_{\text{dry}, 3} \label{eq:return_mass}
\end{align}

Assuming that the fairing is jettisoned during the previous phases of the flight, 
a constraint on the maximum heat flux that the payload undergoes must be included in the formulation from Phase~\StageThree{} to Phase~\StageFourTwo{}.
Specifically, 
modeling the heat flux as a free molecular flow acting on a plane surface perpendicular to the relative velocity \cite{vega2014manual}, the constraint is expressed as
\begin{equation}
    \dot{Q} = \frac{1}{2} \rho v_{\text{rel}}^3 \leq \dot{Q}_{\text{max}}
    \label{eq:heat_flux_nonlinear}
\end{equation}

The goal of the optimization is to minimize the propellant consumed during the fourth stage flight.
Indeed, if the guidance algorithm minimizes the propellant mass consumption, not only the system will feature better responsiveness in off-nominal conditions, but it will also be able to carry additional payload mass, as it requires smaller fuel tanks.
The objective can be equivalently formulated as maximizing the final mass. Thus, the OCP that must be solved to retrieve the optimal trajectory and control in the unperturbed scenario is
\begin{subequations} \label{eq:P_nonlinear_det}
    \begin{align}
        \mathcal{P}_{A} : \;\;  
        \min_{\bm{x}, \bm{\hat{T}}, \Delta t} 
        &\quad -m(\tf)
        \label{eq:objective_m_minimize}
        \\
        \text{s.t.}
        &\quad \text{%
        \eqref{eq:thrust_direction_equality_path_con}--\eqref{eq:heat_flux_nonlinear}}
    \end{align}
\end{subequations}

\section{Convex Formulation of the Deterministic OCP}
\label{sec:cvx}

This section outlines the convexification strategy devised to convert the deterministic OCP in Eq.~\eqref{eq:P_nonlinear_det} into a sequence of convex problems that quickly converges to the optimal solution.


\subsection{Change of Variables}

The coupling of the control variable $\bm{\hat{T}}$ with the state variable $m$ in Eqs.~\eqref{eq:original_ODE_r}--\eqref{eq:original_ODE_m} may lead to numerical issues, such as high-frequency jitters in the control law \cite{liu2015entry}.
Therefore by introducing a new control variable
\begin{equation}
    \bm{u} = \frac{T}{m} \bm{\hat{T}} \label{eq:control_vector_cvx}
\end{equation}
the dynamics are control-affine, that is, they can be written as
\begin{equation}
    \bm{f} = \tilde{\bm{f}}(\bm{x}, t) + \tilde{B} \bm{u}
    \label{eq:ODE_split}
\end{equation}

By replacing Eq.~\eqref{eq:control_vector_cvx} in Eqs.~\eqref{eq:original_ODE_r}--\eqref{eq:original_ODE_m}, the equations of motion become
\begin{align}
    \dot{\bm{r}} &= \bm{v} \label{eq:affine_ODE_r} \\
    \dot{\bm{v}} &= -\frac{\mu}{r^3} \bm{r} + \bm{u} 
    - \frac{D}{m} \bm{\hat{v}}_{\text{rel}} 
    \label{eq:affine_ODE_v} \\
    \dot{m} &= -\dot{m}_e \label{eq:affine_ODE_m}
\end{align}
and the $\tilde{B}$ matrix in Eq.~\eqref{eq:ODE_split} is
$\tilde{B} = [\bm{0}_{3 \times 3} \;
    \bm{I}_{3 \times 3} \;
    \bm{0}_{1 \times 3}]^T$
with $\bm{0}_{m \times n}$ and $\bm{I}_{m \times n}$ denoting the null and identity matrix of size $m \times n$.

The new control variable must satisfy Eq.~\eqref{eq:thrust_direction_equality_path_con}, which is now reformulated as
\begin{equation}
    u_x^2 + u_y^2 + u_z^2 = u_N^2
	\label{eq:thrust_direction_equality_path_con_new}
\end{equation}
where $u_N$ is an additional variable defined as
\begin{equation}
    u_N = \frac{T}{m} \label{eq:u_N}
\end{equation}
In practice, Eq.~\eqref{eq:u_N} must be enforced as a nonlinear path constraint.
This is the price to pay to obtain control-affine dynamics, but the advantages of this change of variables outweigh the disadvantages, as the complexity of the problem is reduced compared to a coupling of state and control variables in the dynamics.

\subsection{Constraint Relaxation}
\label{subsec:constraint_relaxation}
Equation~\eqref{eq:thrust_direction_equality_path_con_new} is a nonlinear equality constraint that must be convexified somehow.
This constraint belongs to a class of constraints that is suitable for a lossless relaxation \cite{acikmese2007convex,accikmecse2011lossless}.
Specifically, by replacing the equality sign with an inequality sign, a second-order cone constraint is obtained
\begin{equation}
    u_x^2 + u_y^2 + u_z^2 \leq u_N^2
	\label{eq:thrust_direction_cone_con}
\end{equation}
This relaxation is lossless because, despite it defines a larger feasible set, the solution of the convex problem is the same as the original and Eq.~\eqref{eq:thrust_direction_cone_con} is satisfied with the equality sign.
A theoretical proof of the lossless property of this relaxation can be found in Proposition~1 of Ref.~\cite{benedikter2022convex}. 

\subsection{Successive Linearization}
The remaining nonlinear expressions are replaced with first-order Taylor series expansions around a reference solution that is updated at every iteration.

\subsubsection{Equations of Motion}
To account for free-time phases in the optimization, 
we introduce an independent variable transformation from physical time $t$ to a new variable $\tau$ defined, for each phase, over a fixed unitary domain $[0, 1]$.
The relationship between the two independent variables is
\begin{equation}
    t^{(i)} = t_0^{(i)} + (t_f^{(i)} - t_0^{(i)}) \, \tau
\end{equation}
Note that the time dilation $s$ between $t$ and $\tau$ corresponds to the time-length of the phase
\begin{equation}
    s^{(i)} = \frac{d t}{d \tau}  = t_f^{(i)} - t_0^{(i)}
    \label{eq:time_dilation}    
\end{equation}
The time-length $s^{(i)}$ is then introduced as an additional optimization variable for every free-time phase.

The equations of motion \eqref{eq:affine_ODE_r}--\eqref{eq:affine_ODE_m} are control-affine but still nonlinear.
Thus, they are expressed in terms of $\tau$ and linearized around a reference solution $\{\bar{\bm{x}}, \bar{\bm{u}}, \bar{s}\}$
\begin{equation}
    \bm{x}' \vcentcolon= \frac{d \bm{x}}{d \tau} = s \bm{f}\left( \bm{x}, \bm{u}, \tau \right) \approx A \bm{x} + B \bm{u} + \Sigma s + \bm{c}
    \label{eq:linear_ODEs}
\end{equation}
where the following matrices and vectors were introduced
\begin{align}
    A &= \bar{s} \frac{\partial \bm{f}}{\partial \bm{x}} ( \bar{\bm{x}}, \bar{\bm{u}}, \tau ) \label{eq:A_matrix_definition} \\
    B &= \bar{s} \frac{\partial \bm{f}}{\partial \bm{u}} ( \bar{\bm{x}}, \bar{\bm{u}}, \tau ) \label{eq:B_matrix_definition} \\
    \Sigma &= \bm{f} ( \bar{\bm{x}}, \bar{\bm{u}}, \tau ) \label{eq:P_matrix_definition} \\
    \bm{c} &= - (A \bar{\bm{x}} + B \bar{\bm{u}}) \label{eq:C_vector_definition}
\end{align}

An undesired phenomenon associated with linearization is artificial infeasibility, that is, the linearized problem may be infeasible.
Therefore, a virtual control signal $\bm{q}$ is added to the dynamics to guarantee the feasibility of the convex problem, that is
\begin{equation}
    \bm{x}' = A \bm{x} + B \bm{u} + \Sigma s + \bm{c} + \bm{q}
    \label{eq:linear_ODEs_vc}
\end{equation}
Since virtual controls are unphysical variables, their use must be highly penalized by adding the following penalty term to the cost function
\begin{equation}
    J_q = \lambda_q P(\bm{q})
    \label{eq:penalty_vc}
\end{equation}
where $\lambda_q$ is a sufficiently large penalty weight and $P(\bm{q})$ a suitable penalty function, which will be defined after discretization in Section~\ref{subsec:discretization}.

\subsubsection{Boundary Constraints}
The final conditions at the payload release \eqref{eq:final_radius_nonlinear}--\eqref{eq:final_r_reentry_nonlinear} are linearized as
\begin{align}
    \bar{\bm{r}}(t_f) \cdot \bar{\bm{r}}(t_f) 
    + 2 \bar{\bm{r}}(t_f) \cdot (\bm{r}(t_f) - \bar{\bm{r}}(t_f)) 
    &= a_{\text{des}}^2 
    \label{eq:final_radius_linearized} \\
    \bar{\bm{v}}(t_f) \cdot \bar{\bm{v}}(t_f) 
    + 2 \bar{\bm{v}}(t_f) \cdot (\bm{v}(t_f) - \bar{\bm{v}}(t_f)) 
    &= v_{\text{des}}^2 
    \label{eq:final_velocity_linearized} \\
    \bar{\bm{r}}(t_f) \cdot \bar{\bm{v}}(t_f) 
    + \bar{\bm{v}}(t_f) \cdot (\bm{r}(t_f) - \bar{\bm{r}}(t_f))&
    \nonumber \\
    + \bar{\bm{r}}(t_f) \cdot (\bm{v}(t_f) - \bar{\bm{v}}(t_f)) &= 0
    \label{eq:final_radial_velocity_linearized} \\
    \bar{v}_y(t_f) (x(t_f) - \bar{x}(t_f))
    - \bar{v}_x(t_f) (y(t_f) - &\bar{y}(t_f))
    \nonumber \\
    - \bar{y}(t_f) v_x(t_f) + \bar{x}(t_f) v_y(t_f) &= h_{z, \text{des}}
    \label{eq:final_ainc_linearized_circular}
\end{align}
Likewise, the terminal radius constraint at the splash-down, Eq.~\eqref{eq:final_r_reentry_nonlinear}, is formulated as
\begin{equation}
    \bar{\bm{r}}(t_f^{(\Return{})}) \cdot \bar{\bm{r}}(t_f^{(\Return{})}) 
    + 2 \bar{\bm{r}}(t_f^{(\Return{})}) \cdot (\bm{r}(t_f^{(\Return{})}) - \bar{\bm{r}}(t_f^{(\Return{})})) 
    = R_E^2
    \label{eq:final_r_reentry_linearized}
\end{equation}

Analogously to the dynamics, also the linearization of these constraints may cause artificial infeasibility. 
Thus, we introduce virtual buffer zones.
In particular, Eqs.~\eqref{eq:final_radius_linearized}--\eqref{eq:final_r_reentry_linearized} are grouped into a vector $\bm{\chi} = \bm{0}$ that is then relaxed as $\bm{\chi} = \bm{w}$.
The vector $\bm{w}$ holds all the virtual buffers, which are highly penalized by adding the following penalty term to the cost function
\begin{equation}
    J_w = \lambda_w \norm{\bm{w}}_1
\end{equation}
with $\lambda_w$ denoting the penalty weight of the virtual buffers. 


\subsubsection{Path Constraints}
The heat flux constraint \eqref{eq:heat_flux_nonlinear} is another nonlinear expression that must be linearized.
By considering its Taylor series expansion, we obtain
\begin{equation}
    \dot{\bar{Q}} + \frac{\partial \dot{\bar{Q}}}{\partial \bm{r}} \cdot (\bm{r} - \bar{\bm{r}}) + \frac{\partial \dot{\bar{Q}}}{\partial \bm{v}} \cdot (\bm{v} - \bar{\bm{v}}) \leq \dot{Q}_{\text{max}}
    \label{eq:heat_flux_linearized}
\end{equation}
where
\begin{align}
    \frac{\partial\dot{\bar{Q}}}{\partial \bm{r}} &= \frac{1}{2} \frac{d \bar{\rho}}{d \bm{r}} \bar{v}_{\text{rel}}^3 + \frac{3}{2} \bar{\rho} \bar{v}_{\text{rel}} \bm{\omega}_E \times \bm{\bar{v}}_{\text{rel}} \label{eq:dQ_dr} \\
    \frac{\partial\dot{\bar{Q}}}{\partial \bm{v}} &= \frac{3}{2} \bar{\rho} \bar{v}_{\text{rel}} \bm{\bar{v}}_{\text{rel}} \label{eq:dQ_dv}
\end{align}

Likewise, the auxiliary variable $u_N$ must satisfy Eq.~\eqref{eq:u_N} at every time.
This is a nonlinear path constraint and it is linearized as
\begin{equation}
    u_N = \frac{T_{vac} - \bar{p} A_e}{\bar{m}} \left(2 - \frac{m}{\bar{m}} \right) - 
    \frac{A_e}{\bar{m}} \frac{d \bar{p}}{d \bm{r}} \cdot (\bm{r} - \bar{\bm{r}})     
    \label{eq:u_N_path_con_linearized}
\end{equation}

\subsection{Trust Region on Time-Lengths}

To enhance convergence, a soft trust region constraint is imposed on the time-lengths $s$ of Phases \CoastingFour{} and \StageFourTwo{}, which resulted particularly sensitive to diverging phenomena \cite{benedikter2022convex}.
The trust region constraint is
\begin{equation}
    | s^{(i)} - \bar{s}^{(i)} | \leq \delta^{(i)} \qquad i = \CoastingFour{}, \StageFourTwo{}
    \label{eq:Dt_trust_region}
\end{equation}
The trust radii $\delta^{(i)}$ are additional optimization variables bounded in a fixed interval $[0, \smash{\delta_{\text{max}}^{(i)}}]$ and slightly penalized in the cost function via the following penalty terms
\begin{equation}
    J_\delta^{(i)} = \lambda_{\delta}^{(i)} \delta^{(i)}
    \qquad i = \CoastingFour{}, \StageFourTwo{}
    \label{eq:delta_penalty_term}
\end{equation}
A reasonable choice of the upper bound \smash{$\delta_{\text{max}}^{(i)}$} is between 1\% and 10\% of $\bar{s}^{(i)}$.
Instead, the penalty weights $\lambda_{\delta}$ should be picked as small as possible (e.g., in the range $\smash{10^{-6}} \div \smash{10^{-3}}$) to prevent converging to a suboptimal solution.

The resulting second-order cone programming problem is
\begin{subequations} \label{eq:P_cvx_det}
    \begin{align}
        \mathcal{P}_{B} : \;\;  
        \min_{\substack{\bm{x}, \bm{u}, s, \\ \bm{q}, \delta^{(\CoastingFour{})}, \delta^{(\StageFourTwo{})}}} 
        &\quad -m(\tf) + J_q + J_\delta^{(\CoastingFour{})} + J_\delta^{(\StageFourTwo{})}
        \label{eq:objective_cvx}
        \\
        \text{s.t.}
        &\quad \text{%
        \eqref{eq:x0}, \eqref{eq:final_LAT_reentry}--\eqref{eq:return_mass}, \eqref{eq:thrust_direction_cone_con}, \eqref{eq:linear_ODEs_vc},} \nonumber 
        \\
        &\quad \text{%
        \eqref{eq:final_radius_linearized}--\eqref{eq:final_r_reentry_linearized}, \eqref{eq:heat_flux_linearized}, \eqref{eq:u_N_path_con_linearized}, \eqref{eq:Dt_trust_region}}
    \end{align}
\end{subequations}

\subsection{Discretization}
\label{subsec:discretization}

Problem $\mathcal{P}_{B}$ in Eq.~\eqref{eq:P_cvx_det} is infinite-dimensional, as state and controls are continuous-time functions.
Therefore, to solve numerically the OCP, we employ a $hp$ pseudospectral discretization method to convert it into a finite set of variables and constraints.
A $hp$ transcription combines the advantages of $h$ and $p$ schemes, as it can introduce mesh nodes near potential discontinuities and leverage exponential convergence rate in each subinterval, resulting in a computationally efficient and accurate discretization \cite{darby2011hp}.

Specifically, in a $hp$ scheme, the time domain of each phase is split into $h$ subintervals, or \emph{segments}, and the differential constraints are enforced in each segment via local orthogonal collocation. 
The order $p$ of the collocation can vary among the segments, but, for the sake of simplicity, we adopt the same order over all the segments of a phase.
The Radau pseudospectral method (RPM) \cite{garg2011advances} is used as the collocation scheme since it is one of the most accurate pseudospectral methods and avoids redundant control variables at the segment boundaries \cite{garg2010unified}.

The state and control are discretized over the grid, generating a finite set of variables $(\bm{x}_j^k, \bm{u}_j^k)$.
The superscript $k$ denotes the $k$-th segment, while the subscript $j$ refers to the $j$-th node of the segment.
In each segment $[\tau_{k}, \tau_{k+1}]$, the state and control signals are approximated via Lagrange polynomial interpolation and the derivative of the state approximation is imposed to be equal to the right-hand side of the equations of motion \eqref{eq:linear_ODEs_vc} at the LGR collocation points
\begin{equation}
    \sum_{j = 1}^{p + 1} D_{ij}^k \bm{x}_j^k = \frac{\tau_{k+1} - \tau_{k}}{2} \bm{f}_i^k 
    \qquad 
    i = 1, \dots, p 
    \label{eq:collocation_radau}
\end{equation}
where $\bm{f}_i^k = A_i^k \bm{x}_i^k + B_i^k \bm{u}_i^k + \Sigma_i^k s + \bm{c}_i^k + \bm{q}_i^k$.
The LGR differentiation matrix $D^k$ can be efficiently computed via barycentric Lagrange interpolation \cite{berrut2004barycentric}.

Note that also the virtual control $\bm{q}$ and the linearization matrices \eqref{eq:A_matrix_definition}--\eqref{eq:C_vector_definition} are discretized over the mesh.
The penalty function $P(\bm{q})$ in Eq.~\eqref{eq:penalty_vc} is chosen as the 1-norm of the vector containing all the discrete virtual control variables $\bm{q}_i^k$.

\subsection{Convergence Criteria}

Problem $\mathcal{P}_{B}$ in Eq.~\eqref{eq:P_cvx_det} must be solved sequentially, i.e., recursively updating the reference values $\{\bar{\bm{x}}, \bar{\bm{u}}, \bar{s}\}$ until convergence.

The update of the reference solution is based on \emph{filtering}, which consists of updating the reference solution with a weighted sum of $K$ previous solutions \cite{benedikter2019convexrendezvous, benedikter2022convex}.
Filtering provides further algorithmic robustness to the successive convexification algorithm since it filters out diverging intermediate iterations by averaging out the reference values over multiple solutions.
In practice, the reference solution for the $i$-th problem is computed as
\begin{equation}
	\bar{x}^{(i)} = \sum_{k = 1}^{K} \alpha_k x^{\max \{0, (i - k)\} }
\end{equation}
where $x^{(i)}$ denotes the solution to the $i$-th problem and $\alpha_k$ the corresponding weight.
As in the authors' previous works \cite{benedikter2019convexrendezvous,benedikter2020autonomous, benedikter2022convex},
three previous solutions are used for the reference solution update, weighted as $\alpha_1 = 6/11$, $\alpha_2 = 3/11$, and $\alpha_3 = 2/11$.

Eventually, the sequential algorithm terminates when all the following criteria are met:
\begin{enumerate}[label=(\roman*)]
\item the difference between the computed solution and the reference one converges below an assigned tolerance $\epsilon_{\text{tol}}$
\begin{equation}
    \norm{\bm{x} - \bar{\bm{x}}}_\infty < \epsilon_{\text{tol}}
    \label{eq:successive_cvx_termination_condition}
\end{equation}
\item the computed solution adheres to the nonlinear dynamics within a tolerance $\epsilon_f$
\begin{equation}
    \Big\lVert \sum_{j = 1}^{p + 1} D_{ij}^k \bm{x}_j^k - \frac{\tau_{k+1} - \tau_{k}}{2} \bm{f}_i^k \Big\rVert_\infty < \epsilon_f 
    \label{eq:acceptable_dyn}
\end{equation}
in each phase, for $i = 1, \dots, p$, and $k = 1, \dots, h$; 
\item the virtual buffers of the computed solution are below the dynamics tolerance $\epsilon_{f}$
\begin{equation}
    \norm{\bm{w}}_\infty < \epsilon_{f}
    \label{eq:acceptable_vb}
\end{equation}
\end{enumerate}

\section{Model Predictive Control}
\label{sec:mpc}

Sections~\ref{sec:problem} and \ref{sec:cvx} outlined the optimization problem to solve to find the optimal trajectory and control law in the unperturbed scenario.
However, because of modeling errors or external disturbances, the actual flight path of the launch vehicle inevitably deviates from the nominal solution.
Thus, a closed-loop controller must be implemented to compensate for deviations and ensure meeting all mission requirements even in presence of uncertainties.
MPC represents an appealing solution since it provides robustness to uncertainties by repeatedly solving onboard problem $\mathcal{P}_B$ to update the trajectory and control law based on the encountered conditions.



The considered MPC controller is implemented only during the third stage flight, i.e., until stage burnout, which is the MPC stopping criterion.
Nevertheless, the \emph{prediction} horizon (i.e., the time domain of the OCP) extends until the payload release into orbit and the splash-down of the third stage.
At the end of each control cycle, the OCP is updated with the measurements coming from the navigation system.
The update concerns the initial condition of the OCP, i.e., Eq.~\eqref{eq:x0}, which is updated with the measured state.
The starting reference solution $\{\bar{\bm{x}}, \bar{\bm{u}}, \bar{s}\}$ is replaced with the optimal solution computed at the previous time, removing the portion of flight elapsed between the two control steps.

Since the time-length of the prediction horizon reduces at each update due to the receding-horizon implementation of the MPC algorithm, to save computational resources, the size of the discretization grid of Phase~\StageThree{} is reduced linearly at every step from the value in Table~\ref{tab:mesh} until a minimum value of $p = \num{5}$.

\begin{table}[h]
    \centering
    \caption{Discretization segments and order in each phase}
    \label{tab:mesh}
    \medskip
    \begin{tabular}{l c c c c c c}
    \toprule
    & \multicolumn{6}{c}{\bf Phase} \\
    \cmidrule(lr){2-7}
    & \StageThree{} & \CoastingThree{} & \StageFourOne{} & \CoastingFour{} & \StageFourTwo{} & \Return{} \\
    \midrule
    $h$ & 1 & 1 & 1 & 1 & 1 & 5 \\
    $p$ & 19 & 9 & 19 & 19 & 19 & 20 \\
    Nodes & 20 & 10 & 20 & 20 & 20 & 101 \\
    \bottomrule
    \end{tabular}
\end{table}

To guarantee recursive feasibility, i.e., the formulation of a feasible OCP at every time step, the constraint on maximum thermal flux \eqref{eq:heat_flux_linearized} has to be relaxed.
Indeed, under nominal operating conditions, the virtual controls and buffers would prevent artificial infeasibility.
Yet, due to off-nominal initial conditions or external disturbances in the dynamics, the heat flux may be greater than expected at some points in the simulation and, since the nominal heat flux constraint is active at some points in the trajectory, there is little or no margin to keep it below the threshold value.
Thus, the infeasibility is solved by relaxing the heat flux constraint \eqref{eq:heat_flux_linearized} as
\begin{equation}
    \dot{\bar{Q}} + \frac{\partial \dot{\bar{Q}}}{\partial \bm{r}} \cdot (\bm{r} - \bar{\bm{r}}) + \frac{\partial \dot{\bar{Q}}}{\partial \bm{v}} \cdot (\bm{v} - \bar{\bm{v}}) \leq \dot{Q}_{\text{max}} (1 + \delta_{\dot{Q}})
    \label{eq:heat_flux_linearized_relax}
\end{equation}
where $\delta_{\dot{Q}}$ is a non-negative optimization variable that is penalized in the cost function via the term 
\begin{equation}
    J_{\dot{Q}} = \lambda_{\dot{Q}} \delta_{\dot{Q}}
\end{equation}
In the authors' experience, high values should not be assigned to the penalty weight \smash{$\lambda_{\dot{Q}}$} since the optimization would exploit virtual controls and may never meet the converge criterion \eqref{eq:acceptable_dyn}.
Rather, considering a conservative nominal threshold \smash{$\dot{Q}_{\text{max}}$} and assigning a small penalty to its violation appears as a much more effective strategy.

Finally, it is worthwhile mentioning that, in the present work, the time to solve the optimization problem is supposed null, but the simulation should introduce a delay between the measurement of the updated state and the actuation of the optimal control law.
This delay can introduce significant deviations from the predicted trajectory if the computation time is long.
Nevertheless, thanks to the convexification strategy, the time required to solve the OCP is assumed to be sufficiently brief. 
Preliminary numerical results confirm the validity of this hypothesis.



\section{Optimal Covariance Control}
\label{sec:covariance_control}

MPC indirectly provides robustness to the system thanks to the continuous update of the trajectory and control that compensates for the measured deviations from the nominal path.
However, a sounder approach to designing a robust controller consists in explicitly accounting for uncertainties in the optimization process.

Generic in-flight disturbances can be modeled as additive white Gaussian noise by considering the following (linearized) stochastic dynamics
\begin{equation}
    d\bm{x} = (A(t) \bm{x} + B(t) \bm{u} + \bm{c}(t)) dt + G(t) d\bm{w}
    \label{eq:x_SDE}
\end{equation}
where $\bm{w} \in \real^{n_w}$ is a multidimensional Wiener process, i.e., $\ev{d\bm{w}} = 0$ and $\ev{d\bm{w} d\bm{w}^T} = dt$, and $G \in \real^{n_x \times n_w}$ is a matrix that determines how the external disturbances affect the system.
If the initial state distribution is Gaussian, then the solution to Eq.~\eqref{eq:x_SDE} is a Gaussian process, that is, $\bm{x}(t) \sim \mathcal{N}(\bm{\mu}(t), P(t)) \; \forall t \geq t_0$, where $\bm{\mu} = \ev{\bm{x}}$ is the mean state and $P = \ev{(\bm{x} - \bm{\mu})(\bm{x} - \bm{\mu})^T}$ is the covariance matrix.

As for the upper stage control problem, the objective of covariance control is, starting from an initial state covariance, that is,
\begin{equation}
    P(t_0) = P_0
    \label{eq:P0}
\end{equation}
to minimize the dispersion of the splash-down location of the spent third stage, i.e., the covariance at splash-down $P(t_R)$.

\subsection{Control Parameterization}

Let us consider a zero-order hold (ZOH) control law, that is, $\bm{u}(t) = \bm{u}(t_j) \; \forall t \in [t_j, t_{j + 1}]$.
Under this assumption, we can discretize the state over a grid of $N$ nodes,
\begin{equation}
    t_0 < t_1 < \cdots < t_{N - 1} 
    \label{eq:discrete_mesh}
\end{equation}
The discrete-time stochastic system dynamics are
\begin{equation}
    \bm{x}_{j + 1} = A_j \bm{x}_j + B_j \bm{u}_j + \bm{c}_j + G_j \bm{w}_j 
    \label{eq:x_discrete_prop}
\end{equation}
for $j = 0, \dots, N - 2$. 
The terms $A_j$, $B_j$, $\bm{c}_j$ are
\begin{align}
    A_j &= \Phi(t_j, t_{j + 1}) \\
    B_j &= \int_{t_j}^{t_{j + 1}} \Phi(s, t_{j + 1}) B(s) ds \label{eq:B_j_integral_stochastic} \\
    \bm{c}_j &= \int_{t_j}^{t_{j + 1}} \Phi(s, t_{j + 1}) \bm{c}(s) ds \label{eq:c_j_integral_stochastic}
\end{align}
where $\Phi(t_j, t_{j + 1})$ denotes the state transition matrix from $t_j$ to $t_{j + 1}$,
and $G_j$ is a matrix such that $G_j \bm{w}_j$ is an $n_x$-dimensional random Gaussian vector with covariance
\begin{equation}
    Q_j = \int_{t_j}^{t_{j + 1}} \Phi(s, t_{j + 1}) G(s) G^T(s) \Phi^T(s, t_{j + 1}) ds 
    \label{eq:Q_j_integral_stochastic}
\end{equation}

To actively control the covariance of a linear system, a feedback control law must be implemented.
Specifically, since the SRM is a non-throttleable engine, a multiplicative feedback law is considered, that is
\begin{equation}
    \bm{u} = \left( I + \skewmat{\bm{\alpha}} \right)
    \bm{\bar{\nu}}
    \label{eq:u_cl}
\end{equation}
where the skew matrix $\skewmat{\bm{\alpha}}$ is
\begin{equation}
    \skewmat{\bm{\alpha}} = \begin{bmatrix}
        0 & -\alpha_3 & \alpha_2 \\
        \alpha_3 & 0 & -\alpha_1 \\
        -\alpha_2 & \alpha_1 & 0
    \end{bmatrix}
\end{equation}
and the vector $\bm{\alpha}$ is a linear state feedback
\begin{equation}
    \bm{\alpha} = K (\bm{x} - \bm{\bar{\mu}})
\end{equation}
The vectors $\bm{\bar{\mu}}$ and $\bm{\bar{\nu}}$ denote the nominal state and control, which are the solution of the unperturbed problem $\mathcal{P}_B$.
Note that, if the $k$-th component of the $\bm{\alpha}$ vector, denoted by $\alpha^{(k)}$, satisfies $\abs{\alpha^{(k)}} \ll 1$ for $k = 1, 2, 3$, then the matrix $\left[ I + \skewmat{\bm{\alpha}} \right]$ approximates a rotation matrix and $\norm{\bm{u}} \approx \norm{\bm{\bar{\nu}}}$.
Thus, the feedback action only rotates the nominal control $\bm{\bar{\nu}}$ without altering its magnitude.

\subsection{Covariance Matrix Propagation}

After replacing Eq.~\eqref{eq:u_cl} in Eq.~\eqref{eq:x_discrete_prop}, the time-evolution of the state covariance matrix $P$ is described by the following equation
\begin{align}
    P_{j + 1} &= (A_j - B_j \skewmat{\bm{\bar{\nu}}}_j K_j) P_j (A_j - B_j \skewmat{\bm{\bar{\nu}}}_j K_j)^T 
    \nonumber \\
    &+ Q_j 
    \label{eq:P_discrete_prop_closed_loop}
\end{align}
where the matrix $\skewmat{\bm{\bar{\nu}}}$ is
\begin{equation}
    \skewmat{\bm{\bar{\nu}}} = \begin{bmatrix}
        0 & -\bar{\nu}_z & \bar{\nu}_y \\
        \bar{\nu}_z & 0 & -\bar{\nu}_x \\
        -\bar{\nu}_y & \bar{\nu}_x & 0
    \end{bmatrix}
\end{equation}

Since both $K$ and $P$ are optimization variables of the covariance control problem, Eq.~\eqref{eq:P_discrete_prop_closed_loop} is a nonlinear constraint.
Thus, we introduce a new control variable
\begin{equation}
    Y_j = K_j P_j
\end{equation}
that replaces $K_j$ in Eq.~\eqref{eq:P_discrete_prop_closed_loop} to obtain
\begin{align}
    P_{j + 1} &= A_j P_j A_j^T
    + A_j Y_j^T \skewmat{\bm{\bar{\nu}}}_j B_j^T
    - B_j \skewmat{\bm{\bar{\nu}}}_j Y_j A_j^T 
    \nonumber \\
    &+ B_j \skewmat{\bm{\bar{\nu}}}_j Y_j P_j^{-1} Y_j^T \skewmat{\bm{\bar{\nu}}}_j^T B_j^T 
    + Q_j 
    \label{eq:P_discrete_prop_closed_loop_eq_Y}
\end{align}
Equation~\eqref{eq:P_discrete_prop_closed_loop_eq_Y} is still nonlinear in $P_j$ and $Y_j$.
Nevertheless, by replacing the equality sign with the inequality sign, the following expression
\begin{align}
    P_{j + 1} &\succeq A_j P_j A_j^T
    + A_j Y_j^T \skewmat{\bm{\bar{\nu}}}_j B_j^T
    - B_j \skewmat{\bm{\bar{\nu}}}_j Y_j A_j^T 
    \nonumber \\
    &+ B_j \skewmat{\bm{\bar{\nu}}}_j Y_j P_j^{-1} Y_j^T \skewmat{\bm{\bar{\nu}}}_j^T B_j^T 
    + Q_j 
    \label{eq:P_discrete_prop_closed_loop_ineq_Y}
\end{align}
can be written as a semidefinite cone constraint using Schur complement's lemma \cite{boyd2004convex}
\begin{equation}
    \begin{bmatrix}
        P_j & Y_j^T \skewmat{\bm{\bar{\nu}}}_j^T B_j^T \\
        B_j \skewmat{\bm{\bar{\nu}}}_j Y_j & \Pi_j
    \end{bmatrix} \succeq 0
    \label{eq:P_discrete_prop_closed_loop_ineq_Y_sdp}
\end{equation}
where
\begin{align}
    \Pi_j &= P_{j + 1} 
    - (A_j P_j A_j^T
    + A_j Y_j^T \skewmat{\bm{\bar{\nu}}}_j B_j^T
    \nonumber \\
    &- B_j \skewmat{\bm{\bar{\nu}}}_j Y_j A_j^T + Q_j)
\end{align}
The inequality sign in Eqs.~\eqref{eq:P_discrete_prop_closed_loop_ineq_Y} and \eqref{eq:P_discrete_prop_closed_loop_ineq_Y_sdp} implies that the difference between the left-hand side and the right-hand side is a positive semidefinite matrix.

Replacing Eq.~\eqref{eq:P_discrete_prop_closed_loop_eq_Y} with Eq.~\eqref{eq:P_discrete_prop_closed_loop_ineq_Y} is a relaxation, as the covariance matrices propagated with the inequality sign may have greater eigenvalues (hence, greater state uncertainties) than the ones that would have been obtained through forward propagation with the equality sign.
Nevertheless, since the objective of the optimization is to minimize the final covariance, we expect that the optimal solution satisfies Eq.~\eqref{eq:P_discrete_prop_closed_loop_ineq_Y} with the equality sign (i.e., the relaxation is lossless) like in the authors' previous work \cite{benedikter2022covariance}. Yet, this statement should be verified by numerical experiments or by applying Pontryagin's maximum principle and inspecting the first-order necessary conditions.

No control can be actuated on the system during the return phase.
Thus, the covariance matrix time-evolution reduces to
\begin{equation}
    P_{R, j + 1} = A_{R, j} P_{R, j} A_{R, j}^T
    \label{eq:P_R_prop}
\end{equation}
where the subscript $R$ indicates that these terms are associated with the return (Phase~\Return{}).
The return covariance is propagated over an evenly spaced grid of $N_R$ points.
Note that no random disturbance is added to the return dynamics.
Also, we assume that the covariance matrix is continuous across the boundary between Phases~\StageThree{} and \Return{}, that is,
\begin{equation}
    P_{R, 0} = P_{N - 1}
    \label{eq:P_lkg}
\end{equation}  

\subsection{Control Magnitude Constraint}

To ensure that $\norm{\bm{u}} \approx \norm{\bm{\bar{\nu}}}$, we must impose $\abs{\alpha^{(k)}} \ll 1$ for $k = 1, 2, 3$.
In practice, one can set a maximum angle $\alpha_{\text{max}}$ that cannot be exceeded, that is,
\begin{equation}
    \norm{\bm{\alpha}_j}_\infty \leq \alpha_{\text{max}}
    \label{eq:max_alpha}
\end{equation}
However, Eq.~\eqref{eq:max_alpha} cannot be enforced as a hard constraint since $\bm{\alpha}$ is a function of the state, which is normally distributed.
Indeed, $\bm{\alpha}$ is a normally distributed vector itself, 
$\bm{\alpha}_j \sim \mathcal{N}(\bm{0}, K_j P_j K_j^T)$. 
Therefore, Eq.~\eqref{eq:max_alpha} must be posed as a chance constraint, that is
\begin{equation}
    \pr{\norm{\bm{\alpha}_j}_\infty \leq \alpha_{\text{max}}} \geq p
    \label{eq:max_alpha_chance}
\end{equation}
Equation~\eqref{eq:max_alpha_chance} is numerically intractable, but it can be replaced with a tractable expression if we recall that $\norm{\bm{\alpha}_j}_\infty \leq \norm{\bm{\alpha}_j}_2$ and consider the following sufficient condition of
the $L_2$ norm of a normally distributed vector \cite{ridderhof2020chance}
\begin{equation}
    \gamma^2(p) \rho(K_j P_j K_j^T) \leq \alpha_{\text{max}}^2
    \Rightarrow
    \pr{\norm{\bm{\alpha}_j}_2 \leq \alpha_{\text{max}}} \geq p
    \label{eq:max_alpha_chance_sufficient}
\end{equation}
where $\rho(M)$ denotes the spectral radius of a matrix $M$, $\gamma(p)$ is a coefficient that depends on the desired probability level $p$ as
\begin{equation}
    \gamma(p) = 
    \begin{cases}
        \sqrt{2 \ln \left( 1 / (1 - p) \right)} & \text{if } n_\alpha \leq 2 \\
        \sqrt{2 \ln \left( 1 / (1 - p) \right)} + \sqrt{n_\alpha} & \text{if } n_\alpha > 2
    \end{cases}
    \label{eq:gamma}
\end{equation}
and $n_\alpha$ is the number of components of $\bm{\alpha}$.

The spectral radius of the positive semidefinite matrix $K_j P_j K_j^T = Y_j P_j^{-1} Y_j^T$ corresponds to its maximum eigenvalue, which can be characterized in epigraph form $\lambda_{\text{max}}(Y_j P_j^{-1} Y_j^T) \leq \tau_j$, where $\tau_j$ is an auxiliary variable that satisfies
\begin{equation}
    \tau_j I - Y_j P_j^{-1} Y_j^T \succeq 0
    \label{eq:tau_nonlinear}
\end{equation}
Equation~\eqref{eq:tau_nonlinear} is particularly convenient since it can be written as a semidefinite cone constraint by leveraging Schur complement's lemma, that is,
\begin{equation}
    \begin{bmatrix}
        \tau_j I & Y_j \\
        Y_j^T & P_j
    \end{bmatrix}
    \succeq 0
    \label{eq:tau_sdp_cone}
\end{equation}

Replacing the spectral radius with $\tau_j$, the sufficient condition in Eq.~\eqref{eq:max_alpha_chance_sufficient} can be written as a linear constraint
\begin{equation}
    \tau_j \leq \alpha_{\text{max}}^2 / \gamma^2(p)
    \label{eq:max_alpha_chance_sufficient_linear}
\end{equation}

\subsection{Objective}

The goal of the controller is to minimize the dispersion of the splash-down point.
Since only the position of the splash-down point is relevant, while velocity and mass at impact are not a concern, the objective function is the maximum eigenvalue of the position elements of the covariance matrix at $t_R$, that is, the first 3 rows and columns of $P_{R, N_R - 1}$.
The maximum eigenvalue is characterized in epigraph form $\lambda_{\text{max}}(P_{R, N_R - 1}^{(1:3 \times 1:3)}) \leq \tau_R$ as
\begin{equation}
    \tau_R I - P_{R, N_R - 1}^{(1:3 \times 1:3)} \succeq 0
    \label{eq:tau_max_R}
\end{equation}

To sum up, given a nominal trajectory $\bm{\bar{\mu}}$ and control $\bm{\bar{\nu}}$, the optimal covariance control problem to solve is
\begin{subequations} \label{eq:P_cov_sdp}
    \begin{align}
        \mathcal{P}_{C} : \;\;  
        \min_{\substack{P_j, P_{R, j}, \\ Y_j, \tau_j, \tau_R}} 
        &\quad \tau_R
        \\
        \text{s.t.}
        &\quad \text{%
        \eqref{eq:P0}, \eqref{eq:P_discrete_prop_closed_loop_ineq_Y_sdp}, \eqref{eq:P_R_prop}, \eqref{eq:P_lkg}, \eqref{eq:tau_sdp_cone}--\eqref{eq:tau_max_R}}
    \end{align}
\end{subequations}
$\mathcal{P}_{C}$ is a semidefinite programming problem that can be solved in polynomial time with interior-point algorithms.

\section{Numerical Results}
\label{sec:results}


The solution algorithms were implemented in C++. 
Gurobi's \cite{gurobi} second-order cone programming solver was used to solve $\mathcal{P}_B$, while MOSEK's \cite{mosek} semidefinite programming solver was used for $\mathcal{P}_C$.
All the computations were carried out on a computer equipped with Intel\textregistered{} Core\texttrademark{} i7-9700K CPU @ \SI{3.60}{\giga\Hz}.

\subsection{Case Study}

The investigated scenario is the same as in Ref.~\cite{benedikter2022convex},
that is, an ascent trajectory from an equatorial launch base to a \SI{700}{\kilo\m} circular orbit with $i_{\text{des}} = \SI{90}{\degree}$ and a splash-down of the third stage targeted at a latitude of $\varphi_{R, \text{des}} = \SI{60}{\degree}$.
Also, the same data are used to model the launch vehicle stages, which are summarized in Table~\ref{tab:engine}.
In addition to these data, the dry mass of the third stage, which is the spent stage mass, is $m_{\text{dry, 3}} = \SI{1326.5}{\kilo\g}$. 
The burn time is denoted by $t_b$, and the vacuum thrust and mass flow rate profiles are modeled as linear functions of time.
As for the aerodynamics, the $C_D = \num{0.38}$ and $S_{\text{ref}} = \SI{9.08}{\m\squared}$ in all phases.
The U.S. Standard Atmosphere 1976 model is used to evaluate the air density and pressure as functions of the altitude~\cite{us1976atm}.

\begin{table}[h]
    \caption{Engine data}
    \label{tab:engine}
    \centering
    \begin{tabular}{c c c c}
        \toprule
        \bf Quantity & \bf Stage 3 & \bf Stage 4 & \bf Unit \\ 
        \midrule
        $t_b$ & \num{104.6} & Free & \si{\s} \\
        $T_{\text{vac, 0}}$ & \num{295.97} & \num{2.45} & \si{\kilo\N} \\
        $T_{\text{vac, f}}$ & \num{221.60} & \num{2.45} & \si{\kilo\N} \\
        $\dot{m}_{e, 0}$ & \num{103.27} & \num{0.79} & \si{\kilo\g\per\s} \\
        $\dot{m}_{e, f}$ & \num{77.32} & \num{0.79} & \si{\kilo\g\per\s} \\
        $A_e$ & \num{1.18} & \num{0.07} & \si{\m\squared} \\
        \bottomrule
    \end{tabular}
\end{table}

Phase~\StageThree{} lasts for the entire burn duration of Stage~3, i.e., $\Delta t^{(\StageThree{})} = t_{b, 3}$.
The duration of the coasting arc after the third stage separation (Phase \CoastingThree{}) is fixed to $\Delta t^{(\CoastingThree{})} = \SI{15.4}{\s}$.
The time-lengths of all the following phases are left to be optimized.
The threshold on the bearable heat flux is set to $\dot{Q}_{\text{max}} = \SI{900}{\W\per\square\m}$.
This value of $\dot{Q}_{\text{max}}$ is lower than the value (\SI{1135}{\watt\per\square\m}) reported in the official VEGA manual \cite{vega2014manual}, to allow for the relaxation of the heat flux constraint as in Eq.~\eqref{eq:heat_flux_linearized_relax} with a relatively small value of the penalty coefficient (\smash{$\lambda_{\dot{Q}} = \num{e-2}$}) in the MPC simulations.

\subsection{Nominal Trajectory}

The nominal trajectory was computed by sequentially solving the convex optimization problem $\mathcal{P}_B$ outlined in Eq.~\eqref{eq:P_cvx_det}.
The same parameters as in Ref.~\cite{benedikter2022convex} were used for the successive convexification algorithm.
Specifically, the penalty weights on the time-lengths, $\smash{\lambda_\delta^{(\CoastingFour{})}}$ and $\smash{\lambda_\delta^{(\StageFourTwo{})}}$, were set to $\smash{\num{e-4}}$, while, to highly penalize solutions that actively exploit virtual variables, we set \smash{$\lambda_q = \lambda_w = \num{e4}$}.
The tolerances on the convergence criteria, Eqs.~\eqref{eq:successive_cvx_termination_condition}--\eqref{eq:acceptable_vb}, were prescribed as \smash{$\epsilon_{\text{tol}} = \num{e-4}$} and \smash{$\epsilon_f = \num{e-6}$}

\begin{figure}[h]
	\centering
    \includegraphics[width=\linewidth,trim={4cm 0.5cm 2cm 4.5cm},clip]{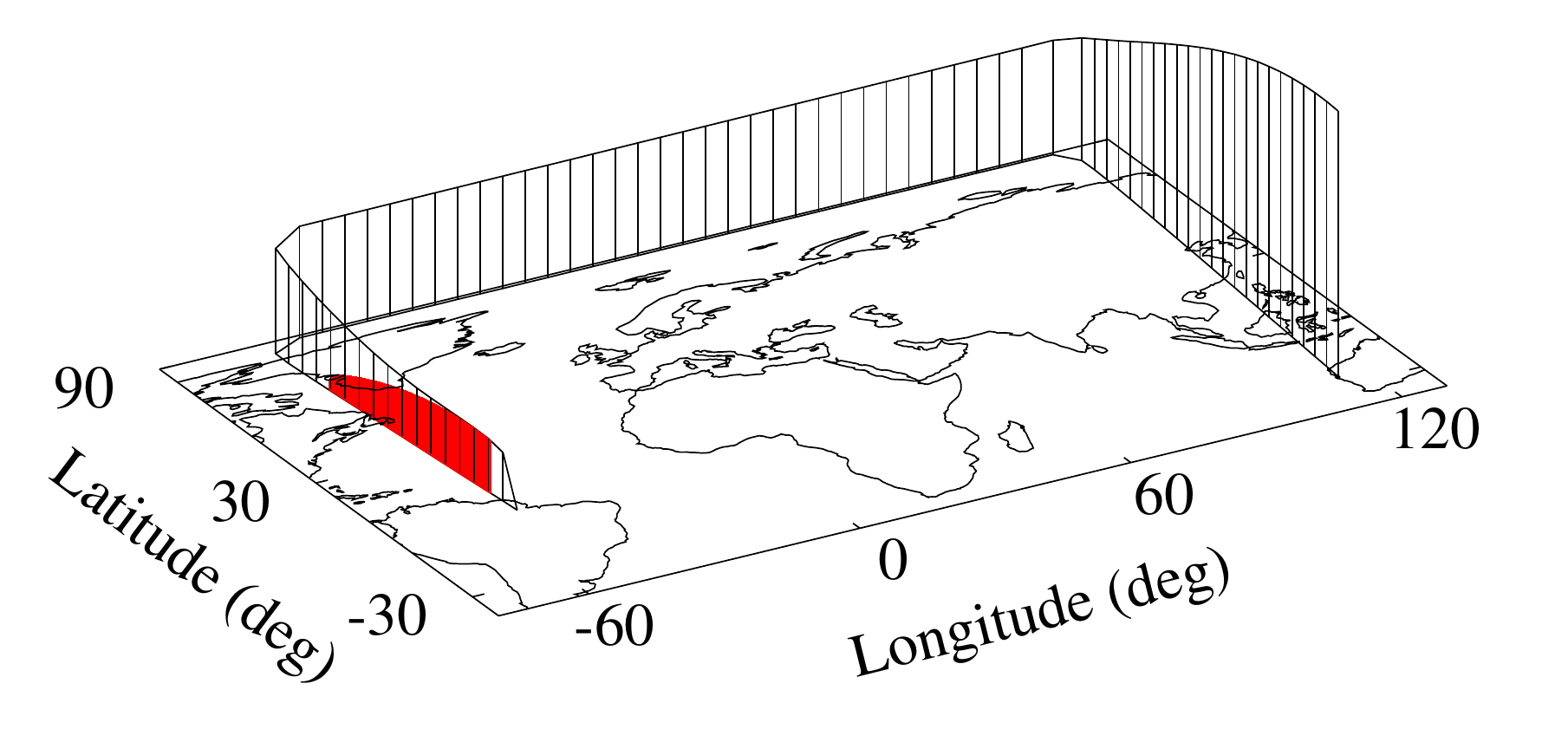}
	\caption{Visualization of the nominal trajectory, with the third stage return flight colored in red.}
	\label{fig:groundtrack_3d}
\end{figure}

First, the entire ascent, from Stage~1 ignition until payload release into orbit was optimized neglecting the return of the spent third stage.
Thanks to the robustness of the convex optimization algorithm, even a rough initial guess leads to the optimal solution in a short time (\SI{25}{s} on average)\footnote{A deeper analysis of the performance and efficiency of the convex optimization of the ascent trajectory is available in Ref.~\cite{benedikter2022convex}}.
Then, the return trajectory was simulated and the splash-down point was gradually moved from the unconstrained latitude to $\varphi_{R, \text{des}}$.
The resulting solution features a nominal payload mass of \SI{1400.1}{\kilo\g} and the corresponding trajectory is illustrated in Fig.~\ref{fig:groundtrack_3d}.

\subsection{Stochastic Control}

In the stochastic optimal control problem, two sources of uncertainty are considered: dispersion of the initial conditions and in-flight random disturbances.
The initial state is normally distributed and features a diagonal covariance matrix.
Specifically, at the initial time, the position variance is \SI{e4}{\square\m}, the velocity variance is \SI{e2}{\square\m\per\square\s}, and the mass variance is \SI{1}{\square\kilo\g}.
On the other hand, the external disturbance is a random perturbing acceleration (hence, the equations of position and mass are unperturbed).
Thus, the dimension of the Wiener process is $n_w = 3$ and the $G$ matrix structure is
\begin{equation}
    G = \begin{bmatrix}
        0_{3 \times 3} &
        g_v I_{3} &
        0_{1 \times 3}^T
    \end{bmatrix}^T
    \label{eq:G_matrix_structure}
\end{equation}
where $g_v$ determines the intensity of the perturbation.
In particular, three levels of in-flight disturbance intensity were considered as reported in Table~\ref{tab:inflight_noise_levels} and called Low~(L), Medium~(M), and High~(H).

\begin{table}[h]
    \caption{In-flight disturbance intensities}
    \label{tab:inflight_noise_levels}
    \centering
    \begin{tabular}{c c c c c}
        \toprule
        \bf Parameter & \bf Case L & \bf Case M & \bf Case H & \bf Unit \\ 
        \midrule
        $g_v$ & \num{0.1} & \num{0.5} & \num{1} & \si{\m\per\s^{3/2}} \\
        \bottomrule
    \end{tabular}
\end{table}

The stochastic control problem was discretized over an evenly spaced grid of $N = 161$ nodes in Phase~\StageThree{} and $N_R = 1001$ nodes in Phase~\Return{}.
The maximum control angle was set to $\alpha_{\text{max}} = \SI{5}{deg}$ and the corresponding chance constraint in Eq.~\eqref{eq:max_alpha_chance} was imposed with a probability level $p = 0.95$.
Problem~$\mathcal{P}_C$ was solved in approximately \SI{40}{\s}.

\subsection{Monte Carlo Campaigns}

A Monte Carlo analysis of the combined effect of off-nominal initial conditions and in-flight disturbance was carried out for both the MPC algorithm and the stochastic controller resulting from problem $\mathcal{P}_C$.
For each disturbance intensity level, 400 independent simulations were carried out.
Each simulation consisted of a numerical integration of the nonlinear stochastic differential equations, i.e., the original Eqs.~\eqref{eq:original_ODE_r}--\eqref{eq:original_ODE_m} perturbed by the Wiener process multiplied by the $G$ matrix in Eq.~\eqref{eq:G_matrix_structure}, using the Euler-Maruyama method.

The update frequency of the MPC was set to \SI{1}{\Hz}, meaning that the OCP is solved every $T = \SI{1}{\s}$.
This is the expected highest achievable rate on a standard flight hardware. 
Note that the solution of the OCP requires approximately \SI{1}{\s} on the employed working station, but significant speed-ups may be attained through custom code optimization or dedicated hardware. 

\begin{figure}[!b]
    \centering
    \begin{minipage}{\columnwidth}
        \centering
        \includegraphics[width=1\linewidth]{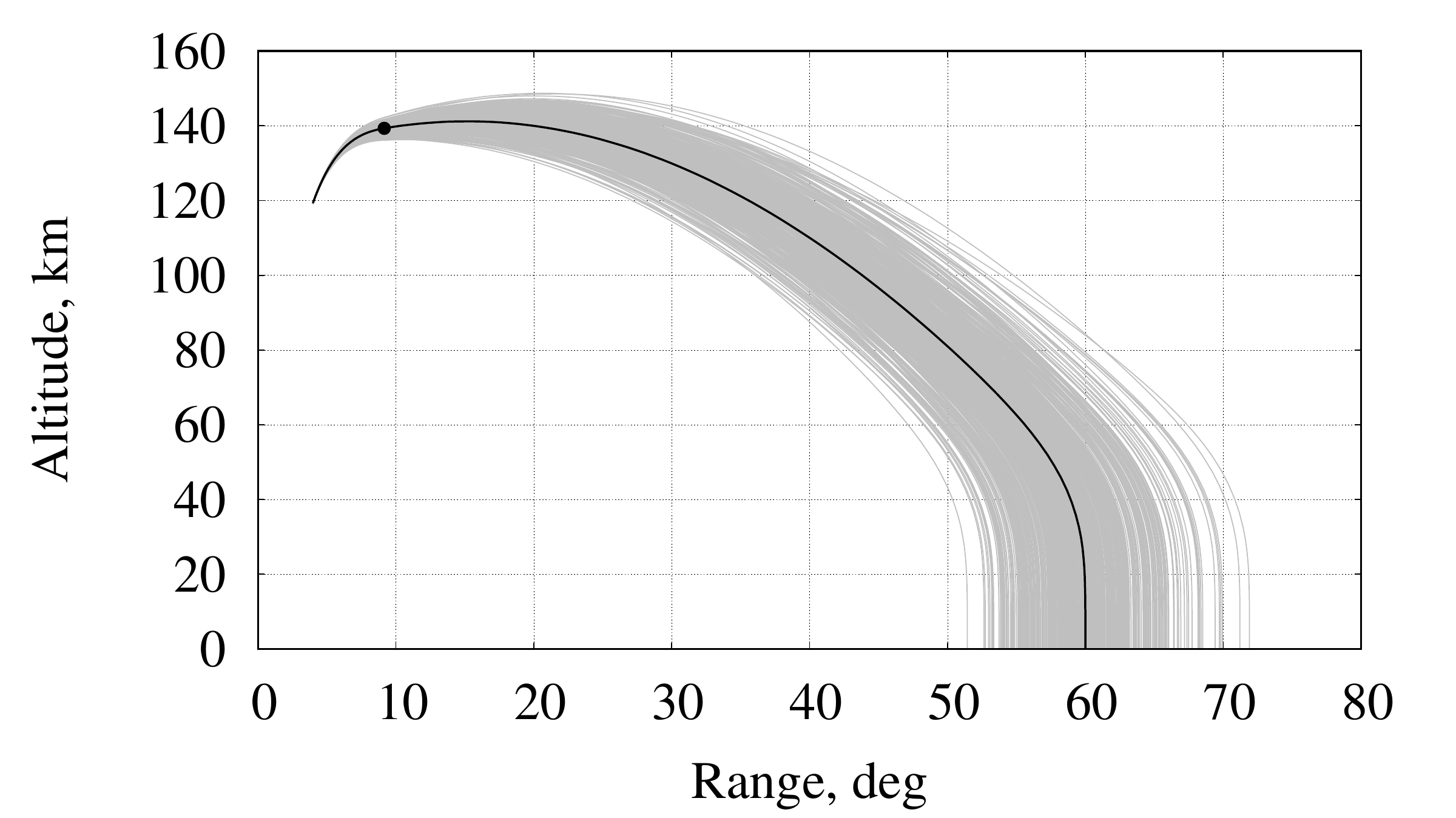}
        \subcaption{Open-loop}
        \label{fig:MC_altitudes_OL}
    \end{minipage}
    
    \centering
    \begin{minipage}{\columnwidth}
        \centering
        \includegraphics[width=1\linewidth]{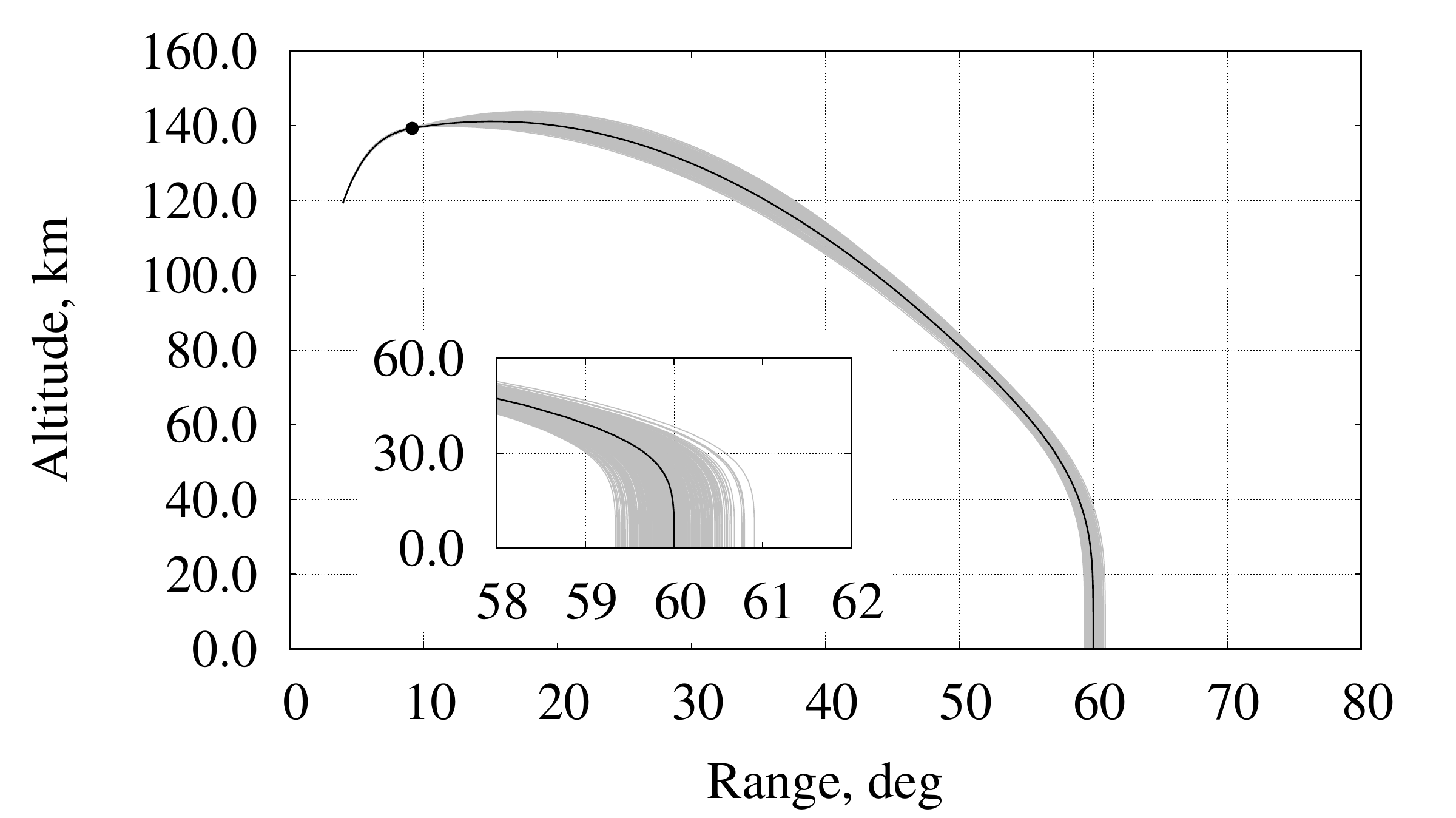}
        \subcaption{Closed-loop (MPC)}
        \label{fig:MC_altitudes_MPC}
    \end{minipage}
    
    \centering
    \begin{minipage}{\columnwidth}
        \centering
        \includegraphics[width=1\linewidth]{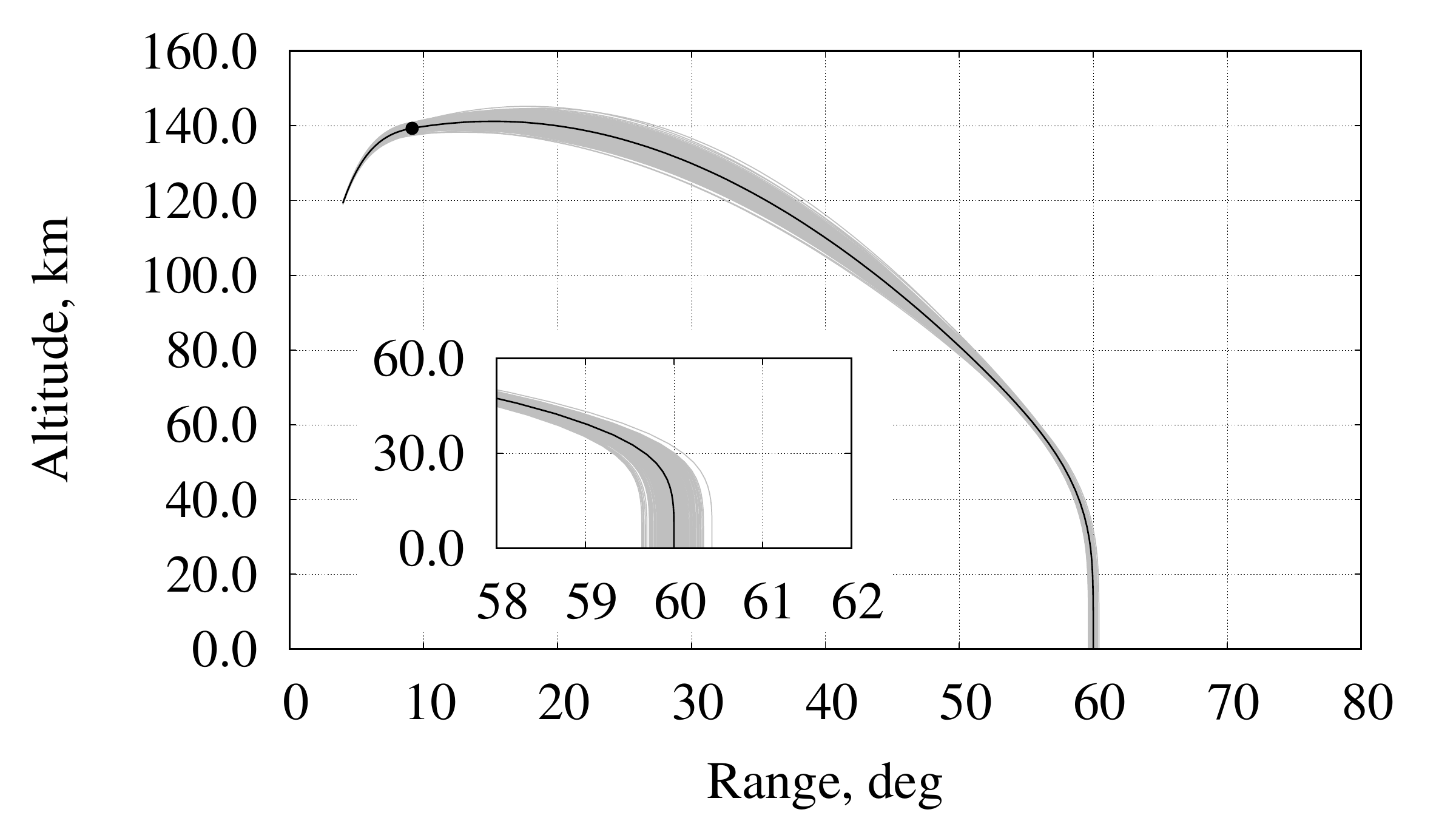}
        \subcaption{Closed-loop (Stochastic control)}
        \label{fig:MC_altitudes_CL}
    \end{minipage}
    
    \caption{Altitude profiles for Case H.}
    \label{fig:MC_altitudes}
\end{figure}

Figure~\ref{fig:MC_altitudes} reports the altitude profile of the third stage from ignition until splash-down versus the range, i.e., the distance from the launch site measured in degrees, in case H.
The black line denotes the nominal flight path, i.e., the solution of Problem~$\mathcal{P}_B$, while the gray lines are the Monte Carlo simulations.
In the Monte Carlo simulations, the Gaussian noise is added to the equations of motion only during the propelled flight of the third stage; after burnout, which is marked by a black dot in the plots, no perturbation is included and the ballistic return is simulated via fourth-order Runge-Kutta integration of the ordinary differential equations.
It is apparent from Fig.~\ref{fig:MC_altitudes_OL} that the dispersion of the splash-down point is huge if no feedback controller is implemented.
Instead, if an MPC algorithm or the stochastic closed-loop controller is implemented, the splash-down dispersion reduces to acceptable values.

\begin{figure}[!b]
    \centering
    \includegraphics[width=0.5\linewidth]{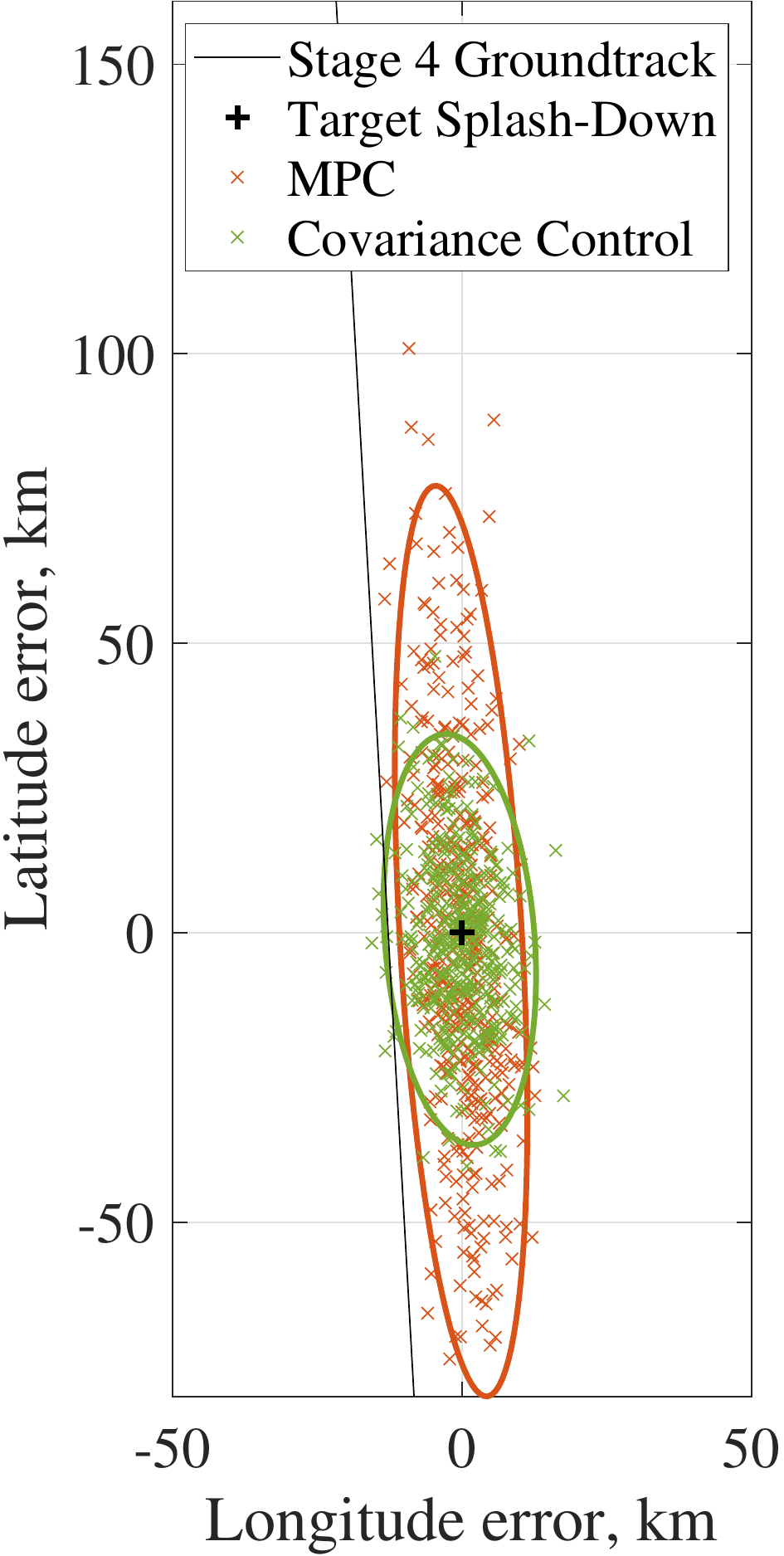}
    \caption{Splash-down footprints for Case H.}
    \label{fig:footprints}
\end{figure}

Figure~\ref{fig:footprints} shows the splash-down footprints attained via MPC and stochastic control.
The footprints are defined as the 95\%-confidence ellipses of the splash-down point.
Stochastic control achieves smaller dispersion than MPC, as the latitude error standard deviation is almost double in the latter case.
Instead, the longitude dispersion is comparable for both control strategies.

\begin{table*}
    \begin{center}
        \caption{Results of the Monte Carlo campaigns}
        \label{tab:MC_results}
        \begin{tabular}{c c c c c c c c c c c c}
        \toprule
        \multirow{2}[2]{*}{\shortstack[c]{\bf Method}} &
        \multirow{2}[2]{*}{\shortstack[c]{\bf Case}} &
        \multicolumn{4}{c}{\bf Payload, kg} &
        \multicolumn{4}{c}{\bf $\varphi_{R}$, deg} &
        \multicolumn{2}{c}{\bf Footprint axes, km}
        \\
        \cmidrule(lr){3-6} 
        \cmidrule(lr){7-10} 
        \cmidrule(lr){11-12} 
        & & Min & Mean & Max & $\sigma$ 
        & Min & Mean & Max & $\sigma$ 
        & Major & Minor 
        \\
        \midrule
        \multirow{3}[0]{*}{\shortstack[c]{MPC}}
        & L &
        \num{1397.7} & \num{1400.0} & \num{1402.2} & \num{0.71} & \num{59.93} & \num{60.00} & \num{60.10} & \num{0.03} & \num{15.09} & \num{2.24} \\
        & M &
        \num{1386.4} & \num{1399.9} & \num{1408.5} & \num{3.48} & \num{59.49} & \num{60.00} & \num{60.44} & \num{0.14} & \num{76.96} & \num{10.64} \\
        & H &
        \num{1375.6} & \num{1399.2} & \num{1418.5} & \num{6.81} & \num{59.34} & \num{59.99} & \num{60.91} & \num{0.29} & \num{157.57} & \num{21.13} \\
        \midrule
        \multirow{3}[0]{*}{\shortstack[c]{Stochastic \\ Control}}
        & L & 
        \num{1377.8} & \num{1399.8} & \num{1418.2} & \num{7.14} & \num{59.91} & \num{60.00} & \num{60.07} & \num{0.03} & \num{15.27} & \num{6.62} \\
        & M & 
        \num{1375.1} & \num{1400.2} & \num{1423.0} & \num{7.46} & \num{59.81} & \num{60.00} & \num{60.20} & \num{0.07} & \num{39.51} & \num{11.56} \\
        & H & 
        \num{1370.4} & \num{1400.3} & \num{1427.2} & \num{9.98} & \num{59.59} & \num{59.98} & \num{60.47} & \num{0.13} & \num{71.63} & \num{27.73} \\
        \bottomrule
        \end{tabular} 
    \end{center}
\end{table*}

Detailed results of the Monte Carlo campaigns are reported in Table~\ref{tab:MC_results}.
As the intensity of the in-flight disturbance increases, larger dispersions on the splash-down point are observed for both control strategies.

\section{Conclusions}
\label{sec:conclusions}

This paper presented a guidance and control strategy for the upper stage of a launch vehicle.
Specifically, a mindful convexification process was outlined to find an optimal nominal trajectory with low computational complexity.
Then, two control synthesis methods were analyzed.
The first one consists in embedding the convex optimization algorithm in an MPC framework, which, thanks to the onboard recomputation of the flight path and control signal, is inherently robust to model mismatches and random inflight disturbances.
The second approach explicitly accounts for uncertainties in the optimization process by considering stochastic system dynamics, perturbed by additive white Gaussian noise, to find a multiplicative feedback controller that guarantees robustness to the considered characterization of the disturbances.
Thanks to an original convexification process, the associated stochastic optimal control problem could be solved in polynomial time.

The performance of the MPC algorithm and the stochastic controller was compared through a Monte Carlo analysis, with a particular focus on the ability to minimize the splash-down footprint of the spent stage.
While MPC proved to be quite robust thanks to the high update frequencies achievable because of the short solution time of the sequential convex optimization algorithm, the stochastic controller turned out to better reject disturbances and attain a smaller dispersion of the splash-down point.
Also, MPC necessitates quite performing onboard computation hardware, while the stochastic controller can be designed preflight and does not require any onboard computation.

Further tests are necessary to prove the suitability of the proposed stochastic control synthesis method for onboard implementation, but the preliminary results are very encouraging.
Indeed, a computationally efficient and systematic design of a nominal trajectory and a robust controller is a crucial technology for launch vehicle guidance and control algorithms. 

\section*{Acknowledgments}
This work was supported by the Agreement n. 2019-4--HH.0 CUP F86C17000080005 
``Technical Assistance on Launch Vehicles and Propulsion'' 
between the Italian Space Agency and the Department of Mechanical and Aerospace Engineering of Sapienza University of Rome.

\bibliographystyle{ieeetr}
\bibliography{references}

\end{document}